# Субгармонические версии

# теоремы Валирона о целых функциях

Б. Н. Хабибуллин[1*], Ф. Б. Хабибуллин[1]

*[1] Башкирский государственный университет*

*Россия, Республика Башкортостан, г. Уфа, 450076, ул. Заки Валиди, 32.*

*\*Email: khabib-bulat@mail.ru*

Оценивается рост канонического интеграла Адамара–Вейерштрасса по мере конечного порядка на комплексной плоскости через типы считающей и усреднённой считающей функции этой меры.

**Ключевые слова:** положительная мера, считающая функция меры, субгармоническая функция, целая функция, порядок роста, функция конечного типа.

**1. Основные обозначения и определения.** Класс всех *положительных борелевских мер,* определённых на борелевских подмножествах комплексной плоскости $\mathbb{C}$ обозначаем как $\text{Meas}(\mathbb{C})$, а класс мер $\mu \in \text{Meas}(\mathbb{C})$ с носителем supp $\mu$, не содержащим точку $0$, – как $\text{Meas}_0(\mathbb{C})$. Для числа $t \in (0, +\infty)$ и меры $\mu \in \text{Meas}(\mathbb{C})$ через $\mu^{\text{rad}}(t)$ обозначаем её *считающую функцию*, равную значениям $\mu(\{z \in \mathbb{C}: |z| \leq t\}) - \mu$-мера замкнутых кругов радиусов $t$ с центром в нуле. Для меры $\mu \in \text{Meas}_0(\mathbb{C})$ будем рассматривать и *усредненную* (проинтегрированную) *считающую функцию*

$$N_\mu(t) := \int_0^t \frac{\mu^{\text{rad}}(t)}{t} \, dt, \qquad t \geq 0. \qquad (1)$$

При $r_0 \in (0, +\infty)$ для функции $f: [r_0, +\infty) \to [0, +\infty)$ величину

$$\text{ord}[f] := \limsup_{r \to +\infty} \frac{\ln^+ f^+(r)}{\ln r} \in [0, +\infty], \quad \text{где} \quad f^+(r) := \max\{f(r), 0\}, \qquad (2)$$

называем *порядком роста* функции $f$ на бесконечности или, для краткости, просто *порядком* $f$. При $\rho \in (0, +\infty)$ для функции $f: [r_0, +\infty) \to [0, +\infty)$ величину

$$\text{type}_\rho[f] := \limsup_{r \to +\infty} \frac{f^+(r)}{r^\rho} \in [0, +\infty] \qquad (3)$$

называем *типом функции $f$ при порядке $\rho$* на бесконечности или просто *типом $f$ при порядке $\rho$*. Мера $\mu \in \text{Meas}(\mathbb{C})$ *порядка $\rho$ и соответственно типа $\sigma$ при порядке $\rho$*, если таковой является её считающая функция $\mu^{\text{rad}}$, т.е. $\text{ord}[\mu] := \text{ord}[\mu^{\text{rad}}] = \rho$ в обозначении (2) и соответственно $\text{type}_\rho[\mu] := \text{type}_\rho[\mu^{\text{rad}}] = \sigma$ в обозначении (3).

*Класс всех субгармонических функций* на $\mathbb{C}$ обозначаем как $\mathrm{sbh}(\mathbb{C})$, а её подкласс, не содержащий функцию, тождественно равную $-\infty$, - через $\mathrm{sbh}_*(\mathbb{C})$. Для функции $u \in \mathrm{sbh}_*(\mathbb{C})$ через $\mu_u \coloneqq \frac{1}{2\pi}\Delta u$, где оператор Лапласа $\Delta$ действует в смысле обобщённых функций, обозначаем *меру Рисса* функции $u$. Функция $u \in \mathrm{sbh}(\mathbb{C})$ порядка $\rho$ и соответственно типа $\sigma$ при порядке $\rho$, если в обозначениях (2) и (3) такова функция $M_u(r) \coloneqq \max\{u(z): |z| = r\}$, т.е. $\mathrm{ord}[u] \coloneqq \mathrm{ord}[M_u] = \rho$ и соответственно $\mathrm{type}_\rho[u] \coloneqq \mathrm{type}_\rho[M_u] = \sigma$ в обозначениях (2) и (3). Целая функция $g \neq 0$ на $\mathbb{C}$ - функция порядка $\rho$ и соответственно типа $\sigma$ при порядке $\rho$, если такова функция $\ln|g| \in \mathrm{sbh}_*(\mathbb{C})$, т.е., с некоторой вольностью в обозначениях, $\mathrm{ord}[g] \coloneqq \mathrm{ord}[\ln|g|] = \rho$, $\mathrm{type}_\rho[g] \coloneqq \mathrm{type}_\rho[\ln|g|] = \sigma$.

**2. Основной результат.** Пусть $\mathbb{N}_0 \coloneqq 0,1,2,\ldots$ - «французское» множество натуральных чисел, $q \in \mathbb{N}_0$. *Первичное ядро Вейерштрасса рода $q$* определяется как [1; 4.1]

$$K_q(z) \coloneqq \ln|1-z| + \sum_{k=1}^{q} \mathrm{Re}\, \frac{z^k}{k}, \text{ где, как обычно, } \sum_{k=1}^{0} \cdots \coloneqq 0. \tag{4}$$

Пусть мера $\mu \in \mathrm{Meas}_0(\mathbb{C})$ конечного типа при порядке $\rho$. Тогда определен *канонический интеграл Адамара–Вейерштрасса по мере $\mu$ рода $q = [\rho]$ — целая часть $\rho$* [1; Лемма 4.4], задаваемый с помощью ядра Вейерштрасса (4) и обозначаемый как

$$U_q^\mu(z) \coloneqq \int_\mathbb{C} K_q\left(\frac{z}{w}\right) \mathrm{d}\mu(w), \quad z \in \mathbb{C}. \tag{5}$$

Здесь мы получаем верхние оценки для канонического интеграла Адамара–Вейерштрасса (5) через рост считающей функции $\mu^{\mathrm{rad}}$ и усредненной считающей функции (1). Первый вариант со считающей функцией нулей (корней) целой функции был получен еще в начале XX в. Ж. Валироном в [3]. В последние годы тонкие и глубокие результаты в терминах уточненного порядка для целых функций в этом направлении были получены А. Ю. Поповым [4] и Ф. С. Мышаковым [5]. Полагаем

$$\mathcal{M}_q(r) \coloneqq \max_{|z|=r} K_q(z), \; r \geq 0. \tag{6}$$

Поскольку ядро Вейерштрасса (4) – субгармоническая функция, а функция (6) непрерывна, то $\mathcal{M}_q(r)$ – *возрастающая выпуклая от $\ln$ положительная непрерывная функция* [2; Теорема 2.6.8]. В частности, при всех $r > 0$ у функции (6) *существуют левые и правые производные, и они совпадают вне не более чем счетного множества точек на* $(0,+\infty)$. Это позволяет использовать в интегралах одну из производных. Для определенности далее используем в интегралах и дифференциалах правую производную, а обозначаем её как обычную производную $\mathcal{M}_q'$. При этом в силу выпуклости от $\ln$ функции (6) функция $0 \leq t \mapsto t\, \mathcal{M}_q'(t)$ − *возрастающая и положительная*. Следуя [4]-[5],

$$S(\rho) := \int_0^{+\infty} r^{-\rho}\, \mathcal{M}_q'(r)\mathrm{d}r = \rho \int_0^{+\infty} r^{-\rho-1}\, \mathcal{M}_q(r)\, \mathrm{d}r, \qquad (7)$$

где величина (7) конечна при нецелом $\rho \notin \mathbb{N}_0$.

**Теорема.** *Пусть* $\mu \in \mathrm{Meas}_0(\mathbb{C})$, $\mathrm{type}_\rho[\mu] < +\infty$, $q := [\rho]$. *Тогда*

**1.** *Справедливы соотношения*

$$U_q^\mu(z) \leq \int_0^{+\infty} \mathcal{M}_q\left(\frac{r}{t}\right) \mathrm{d}\mu^{\mathrm{rad}}(t) = \int_0^{+\infty} \mu^{\mathrm{rad}}\left(\frac{r}{t}\right) \mathcal{M}_q'(t)\mathrm{d}t, \qquad r = |z| \geq 0.$$

**2.** *В обозначениях* (1), (6) *справедливы соотношения*

$$U_q^\mu(z) \leq \int_0^{+\infty} \frac{r}{t} \mathcal{M}_q\left(\frac{r}{t}\right) \mathrm{d}N_{\mu_u}(t) = \int_0^{+\infty} N_{\mu_u}\left(\frac{r}{t}\right) \mathrm{d}(t\mathcal{M}_q'(t)), \qquad r = |z| \geq 0.$$

**3.** *При* $\rho \notin \mathbb{N}_0$ *для любой функции* $u \in \mathrm{sbh}_*(\mathbb{C})$ *с* $\mathrm{type}_\rho[u] < +\infty$ *и мерой Рисса* $\mu_u$ *в обозначениях* (1), (7) *имеют место оценки*

$$\mathrm{type}_\rho[u] \leq S(\rho)\,\mathrm{type}_\rho[\mu_u], \quad \mathrm{type}_\rho[u] \leq \rho S(\rho)\,\mathrm{type}_\rho[N_{\mu_u}].$$

**Комментарий.** Соотношения из п. **1** нашей Теоремы использованы для целых функций в виде рядов по последовательности нулей канонического произведения Адамара–Бореля Ж. Валироном [3], установившим для целых функций первую оценку из п. **3**, и А. Ю. Поповым в [4], существенно уточнившим и развившим эти результаты. Для целых функций результат п. **2** и вторая оценка из п. **3** установлены и значительно развиты Ф. С. Мышаковым в [5] даже для уточнённых порядков. В [4]-[5] показано, что эти результаты точны в классе целых функций и последовательностей точек-нулей, значит они тем более точны в классе субгармонических функций и мер конечного типа при порядке $\rho$. Разнообразные, более или менее явные, формулы для функций (6)-(7) и их производных можно найти в работе А. Данжуа [6] и в [4]-[5]. Дополнительную информацию подобного рода для можно извлечь и из работы одного из авторов [7; § 6] и статьи С. Г. Мерзлякова [8].

*Доказательство Теоремы.* **1.** Из (4)-(6) при $z = re^{i\theta} \in \mathbb{C}, r \geq 0$, следует

$$U_q^\mu(z) \leq \max_{|z|=r} \int_{\mathbb{C}} K_q\left(\frac{re^{i\theta}}{te^{i\varphi}}\right) \mathrm{d}\mu(te^{i\varphi}) \leq \int_0^{+\infty} \mathcal{M}_q\left(\frac{r}{t}\right) \mathrm{d}\mu(te^{i\varphi})$$

$$= \int_0^{+\infty} \mathcal{M}_q\left(\frac{r}{t}\right) \mathrm{d}\mu^{\mathrm{rad}}(t) = -\int_0^{+\infty} \mu^{\mathrm{rad}}(t)\, \mathrm{d}\mathcal{M}_q\left(\frac{r}{t}\right), \qquad (8)$$

где последнее равенство получается интегрированием по частям с учетом ограничения $\text{type}_\rho[\mu] < +\infty$. Замена отношения $r/t$ в последнем интеграле одной переменной завершает доказательство п. **1**.

**2.** В тех же обозначениях продолжим цепочку (8): $U_q^\mu(z) \leq$

$$\int_0^{+\infty} \mu^{\text{rad}}(t)\, \mathcal{M}_q'\!\left(\frac{r}{t}\right)\frac{r\,\mathrm{d}t}{t^2} = \int_0^{+\infty} \mathcal{M}_q'\!\left(\frac{r}{t}\right)\frac{r}{t}\,\mathrm{d}\int_0^t \frac{\mu^{\text{rad}}(s)}{s}\mathrm{d}s = \int_0^{+\infty} \mathcal{M}_q'\!\left(\frac{r}{t}\right)\frac{r}{t}\,\mathrm{d}N_\mu(t),$$

что доказывает неравенство из п. **2**. Заметим, что подынтегральная функция $\mathcal{M}_q'\!\left(\frac{r}{t}\right)\frac{r}{t}$ в последнем интеграле – убывающая по $t$ функция. Вновь интегрируя по частям последний интеграл, получаем продолжение (дифференциал $\mathrm{d}$ действует по $t$)

$$=\int_0^{+\infty} N_\mu(t)\,\mathrm{d}\!\left(\mathcal{M}_q'\!\left(\frac{r}{t}\right)\frac{r}{t}\right) = \int_0^{+\infty} N_{\mu_u}\!\left(\frac{r}{t}\right)\mathrm{d}(t\mathcal{M}_q'(t)),$$ что доказывает п. **2**.

**3.** Пусть $\sigma := \text{type}_\rho[\mu_u]$. Тогда для любого сколь угодно малого числа $\varepsilon > 0$ меру Рисса $\mu_u$ функции $u \in \text{sbh}_*(\mathbb{C})$ конечного типа $\sigma$ при порядке $\rho$ можно представить в виде суммы двух мер $\mu_u = \mu_0 + \mu_\infty$, где мера $\mu_0$ сосредоточена в круге достаточно большого радиуса, т.е. имеет компактный носитель и конечна, а мера $\mu_\infty$ удовлетворяет оценке

$$\mu_\infty(r) \leq (\sigma + \varepsilon)r^\rho \text{ при всех } r \geq 0. \tag{9}$$

Функцию $u \in \text{sbh}_*(\mathbb{C})$ при этом можно представить в виде

$$u = U_0^{\mu_0} + U_q^{\mu_\infty} + H, \tag{10}$$

где $H$ − гармоническая функция и $U_0^{\mu_0} + H(z) = o(|z|^\rho), z \to \infty$. Поэтому оба равенства из п. **3** достаточно установить для канонического интеграла Адамара-Вейерштрасса $U_q^{\mu_\infty}$ по мере $\mu_\infty$, удовлетворяющей оценке (9).

Ввиду оценки п. **1** через последний интеграл в обозначении (7) следует $U_q^{\mu_\infty}(re^{i\theta})$

$$\leq \int_0^{+\infty} (\sigma + \varepsilon)\left(\frac{r}{t}\right)^\rho \mathcal{M}_q'(t)\mathrm{d}t = (\sigma + \varepsilon)r^\rho \int_0^\infty t^{-\rho}\,\mathcal{M}_q'(t)\mathrm{d}t = (\sigma + \varepsilon)r^\rho S(\rho).$$

Поделив обе части на $r^\rho$ и устремляя $r$ к $+\infty$, получаем, в силу произвола в выборе числа $\varepsilon > 0$, получаем первую оценку из п. **3**. Для доказательства второй оценки из п. **3**, используя представление, аналогичное (10), можем считать, как в (9), что $N_{\mu_\infty}(r) \leq$

$(\sigma+\varepsilon)r^\rho$ при всех $r\geq 0$. Ввиду оценки п. 2 через последний интеграл в обозначении (7) следует $U_q^{\mu_\infty}(re^{i\theta}) \leq \int_0^{+\infty} N_{\mu_\infty}\left(\frac{r}{t}\right) \mathrm{d}(t\mathcal{M}_q'(t)) \leq \int_0^{+\infty}(\sigma+\varepsilon)\left(\frac{r}{t}\right)^\rho \mathrm{d}(t\mathcal{M}_q'(t)) = (\sigma+\varepsilon)r^\rho \int_0^{+\infty} \rho t^{-\rho-1}\, t\mathcal{M}_q'(t)\mathrm{d}t = (\sigma+\varepsilon)\rho r^\rho S(\rho),$ где использовано интегрирование по частям. Теорема доказана.

# Subharmonic Versions of Valiron's Theorem on Entire Functions

B. N. Khabibullin[1*], F. B. Khabibullin[1]

[1] *Bashkir State University*
*32 Zaki Validi st., 450074 Ufa, Republic of Bashkortostan, Russia.*

*Email: khabib-bulat@mail.ru.*

We estimate the growth of the canonical integral of Hadamard–Weierstraß of measure of finite order on the complex plane by the type of counting function or average counting function of this measure.

**Keywords:** positive measure, counting function of measure, subharmonic function, entire function, order of the growth, function of finite type.

.